\documentclass{conm-p-l}
\usepackage{amssymb,latexsym, amscd}
  \usepackage[all]{xy}

 \newtheorem{theorem}{Theorem}[subsection]
 \newtheorem{cor}[theorem]{Corollary}
 \newtheorem{lemma}[theorem]{Lemma}
 \newtheorem{proposition}[theorem]{Proposition}
 \theoremstyle{definition}
 \newtheorem{definition}[theorem]{Definition}
 \theoremstyle{definition}
 
 \theoremstyle{remark}
 \newtheorem{rem}[theorem]{Remark}
 \numberwithin{equation}{subsection}

\usepackage{latexsym}
\usepackage{amssymb}

\newcommand{\ben}{\begin{equation}}
\newcommand{\een}{\end{equation}}

\newcommand{\integer}{\ensuremath{{\mathbb Z}}}

\newcommand{\real}{\ensuremath{{\mathbb R}}}
\newcommand{\complex}{\ensuremath{{\mathbb C}}}

\newcommand{\U}[1]{\ensuremath{{\mathrm U( #1 )}}}

\newcommand{\Aa}{{\mathcal A}}
\newcommand{\PP}{{\mathcal P}}

\newcommand{\VV}{{\mathcal V}}
\newcommand{\WW}{{\mathcal W}}

\newcommand{\FF}{{\mathcal F}}
\newcommand{\GG}{{\mathcal G}}

\newcommand{\LL}{\mathcal{L}}

\newcommand{\Xx}{\mathsf{X}}
\newcommand{\Gg}{\mathsf{G}}
\newcommand{\Hh}{\mathsf{H}}
\newcommand{\Mm}{\mathsf{M}}
\newcommand{\Ss}{\mathsf{S}}
\newcommand{\Ll}{\mathsf{L}}
\newcommand{\Ff}{\mathsf{F}}
\newcommand{\Kk}{\mathsf{K}}
\newcommand{\LLL}{\mathrm{L}}
\newcommand{\Loop}{\mathsf{L}}

\newcommand{\twoarrows}{\rightrightarrows}

\newcommand{\twodownarrows}{\downdownarrows }
\newcommand{\timests}{\: {}_{t}  \! \times_{s}}

\begin{document}

\title[Loop Groupoids and Twisted Sectors]{Loop Groupoids, Gerbes, and Twisted Sectors on Orbifolds.}
\author{ Ernesto Lupercio and Bernardo Uribe}

\address{Department of Mathematics, University of Wisconsin at Madison, Madison, WI 53706}
\address{Max-Planck-Institut f\"{u}r Mathematik, Vivatsgasse
7, D-53111 Bonn, Germany Postal Address: PO.Box: 7280, D-53072
Bonn} \email{ lupercio@math.wisc.edu \\ uribe@mpim-bonn.mpg.de}

\begin{abstract}
The purpose of this paper is to introduce the notion of loop
groupoid $\Ll \Gg$ associated to a groupoid $\Gg$. After studying
the general properties of $\Ll \Gg$, we show how this notion
provides a very natural geometric interpretation for the twisted
sectors of an orbifold \cite{Kawasaki}, and for the inner local
systems introduced by Ruan \cite{Ruan} by means of a natural
generalization of the concept holonomy of a gerbe.
\end{abstract}

\maketitle

\section{Introduction}

\subsection{} The study of the free loop space $\LL M$ of a space
$M$ without singularities (say for example an smooth manifold) has
proved very important in geometry, topology, representation
theory and string theory. It is a very natural question that of
generalizing this notion to more general spaces, like an orbifold
or the leaf space of a foliation. The purpose of this paper is to
put forward the notion of the \emph{Loop Groupoid} $\Loop \Gg$
associated to a given groupoid $\Gg$, and to show how this notion
provides a very natural framework to understand geometrically
several notions that have appeared recently in the theory of the
Chen-Ruan cohomology of orbifolds, and also in the work of Freed,
Hopkins and Teleman  on the Verlinde Algebra \cite{Freed}. Some
of the constructions that we propose work just as well in the case
of the leaf space of a foliation.

We should point out that in the case in which the groupoid $\Gg$
represents an orbifold, the coarse moduli space of $\Loop \Gg$ is
precisely the loop space of an orbifold defined by Chen in his
work on orbispaces \cite{Chen}.

\subsection{} Let us briefly describe the contents of this paper. In
section 2 we collect some facts about groupoids and set the
notations used in the rest of the paper. In this paper by the
word groupoid we mean a topological category in which every
morphism has an inverse. In most of the paper we will suppose
that our groupoids are \'etale and proper. In section 3 we define
the loop groupoid associated to a given groupoid. This is
topological category but it is far from being \'etale.  In
section 4 we have an algebraic interlude showing how this theory
looks in the case in which the groupoid $\Gg$ is simply a finite
group. Here we also make some remarks relating this results to the
loop group and compact lie groups, and the work of Freed, Hopkins
and Teleman mentioned before. If our groupoid represents a smooth
manifold $M$ then the loop groupoid represents the free loop
space $\LL M$, this is proved in section 5. In this section we
also explain how this groupoid is useful in understanding how the
holonomy of a gerbe on $M$ induces a line bundle on $\LL M$, and
therefore a one dimensional representation of the loop groupoid
$\Loop \Mm$. In section 6 we consider groupoids representing
orbifolds. We show that the twisted sectors \cite{Ruan} of an
orbifold can be understood as the $S^1$ in variant subgroupoid of
the loop groupoid. This helps us to put all the story together
and show how an inner local system in the sense of Ruan
\cite{Ruan} is simply the natural generalization of the holonomy
line bundle of a gerbe that we have defined before in
\cite{LupercioUribe}.

In this paper we try to be reasonably self-contained, but we will
use at several points the results of \cite{LupercioUribe}, and
this paper can be thought of as a continuation of that one. We
strongly recommend \cite{Moerdijk2} for an introduction to the
subject, including the original references. We refer the reader to
\cite{Moerdijk2, LupercioUribe} for the motivations and the
details of those results.

We would like to thank A. Adem, I. Moerdijk, T. Nevins,  M.
Poddar, J. Robbin and Y. Ruan for useful conversations. We would
like to thank specially the referee for many important remarks
and some references. We would also like to thank A. Adem,
J.~Morava and Y. Ruan for organizing this excellent workshop.

\section{Groupoids}

In this section we will briefly summarize the definition of
groupoids as well as gerbes over groupoids; for a more detailed
exposition we recommend the reader to see \cite{Moerdijk1,
CrainicMoerdijk, LupercioUribe}.

A groupoid $\Gg$ in great generality is a category in which every
morphism is invertible. Here we will only consider small
categories.\\ We will denote by $\Gg_0$ and $\Gg_1$ the set of
objects and morphism respectively, and the structure maps by:
       $$\xymatrix{
         \Gg_1 \timests \Gg_1 \ar[r]^{m} & \Gg_1 \ar[r]^i &
         \Gg_1 \ar@<.5ex>[r]^s \ar@<-.5ex>[r]_t & \Gg_0 \ar[r]^e & \Gg_1
         }$$
where $s$ and $t$ are the source and the target maps of
morphisms, $m$ is the composition of two of them whenever the
target of the first equals the source of the second, $i$ gives us
the inverse morphism and $e$ assigns the identity arrow to every
object. We will also call arrows the morphisms of the groupoid.

The groupoid will be called {\it topological (smooth)} if the
sets $\Gg_1$ and $\Gg_0$
 and the structure maps belong to the category of topological spaces (smooth manifolds).
In the case of a smooth groupoid we will also require that the
maps $s$ and $t$ must be submersions, so that $\Gg_1 \timests
\Gg_1 $ is also a manifold.

A topological (smooth) groupoid is called {\it \'{e}tale} if the
source and target maps $s$ and $t$ are local homeomorphisms (local
diffeomorphisms). For an \'{e}tale groupoid we will mean a
topological \'{e}tale groupoid.  We will always denote groupoids
by letters of the type $\Gg,\Hh,\Ss$. We will work with
topological groupoids but the reader can think of them as smooth
if he prefers.

   Due to a fundamental result of Moerdijk and Pronk \cite{MoerdijkPronk1, Pronk}
we know that orbifolds are related to a special kind of \'{e}tale
groupoids, they have the peculiarity that the anchor map $(s,t):
\Gg_1 \to \Gg_0\times \Gg_0$ is proper, groupoids with this
property are called {\it proper}.  The theorem of Moerdijk and
Pronk states that the category of orbifolds is equivalent to a
quotient category of the category of proper \'{e}tale groupoids
after inverting Morita equivalence. We explain this in more
detail now. Whenever we write orbifold, we will choose a proper
\'{e}tale smooth groupoid representing it (up to Morita
equivalence.)

A morphism of groupoids $\Psi: \Hh \to \Gg$ is a pair of maps
$\Psi_i: \Hh_i \to \Gg_i$ $i=0,1$ such that they commute with the
structure maps. The maps $\Psi_i$ will be continuous (smooth)
depending on which category we are working on.

The morphism $\Psi$ is called {\it Morita} if the following
square is a cartesian square .
\begin{eqnarray}
  \xymatrix{
          \Hh_1 \ar[r]^{\Psi_1} \ar[d]_{(s,t)} & \Gg_1 \ar[d]^{(s,t)} \\
          \Hh_0 \times \Hh_0 \ar[r]^{\Psi_0 \times \Psi_0} & \Gg_0 \times \Gg_0
          } \label{Moritasquare}
\end{eqnarray}

Two (abstract) groupoids $\Gg$ and $\Hh$ are Morita equivalent if
there exist another groupoid $\Kk$ with Morita morphisms $\Gg
\stackrel{\simeq}{\leftarrow} \Kk \stackrel{\simeq}{\to} \Hh$. We
will write $B \Gg$ to denote the classifying space of the
groupoid $\Gg$ \cite{Segal}. If two groupoids $\Gg$ and $\Hh$ are
Morita equivalent we write $\Gg \simeq \Hh$, in this case the
classifying spaces are homotopy equivalent $B \Gg \simeq B \Hh$.

On working with topological groupoids there is a further
requirement for a morphism to be called Morita, namely given the
morphism $\Phi : \Hh \to \Gg$ we require the map $s \pi_2:\Hh_0
\: {}_{\Psi_0} \! \times_{t} \Gg_1 \to \Gg_0$ to be an open
surjection. If this is the case for a morphism $\Phi$ between
topological \'{e}tale groupoids, then we necessarily have that
the maps $\Phi_1$ and $\Phi_0$ are both \'{e}tale (see
\cite{Pronk}). As this paper only deals with topological
groupoids, whenever we mention a Morita morphism we will require
this additional condition.

If $G$ is a group, we will write $\overline{G}$ for the groupoid
$* \times G \twoarrows *$ where $m(g,h)=gh$. All the group
actions will be carried out from the right. A morphism from $\Hh$
to $\overline{G}$ determines a unique principal $G$ bundle over a
groupoid $\Hh$ (cf. 4.1.5 \cite{LupercioUribe}).

A gerbe $(\LL,\theta)$ over a \'{e}tale groupoid $\Gg$ (in most
cases it will be an orbifold) is a complex line bundle $\LL$ over
$\Gg_1$ (not over $\Gg$) satisfying the following conditions:
\begin{itemize}
\item $i^*\LL \cong \LL$
\item $\pi_1^*\LL \otimes \pi_2^*\LL \otimes m^*i^* \LL \stackrel{\theta}{\cong}1$
\item $\theta : \Gg_1 \timests \Gg_1 \to  \U{1}$ is a 2-cocycle
\end{itemize}
where $\pi_i$'s are projections and $\theta$ is a trivialization
of the line bundle. The cohomology class $\langle \LL \rangle \in
H^2(\Gg,\underline{\complex^*})$ of $\theta$ is called the
characteristic class of the gerbe ($ \underline{\complex^*}$ is
the sheaf of $\complex^*$ valued functions over the groupoid
$\Gg$).

A {\it connection} on a gerbe $(\LL,\theta)$ over an \'{e}tale
groupoid $\Gg$ consist of a real valued 0-form $ \theta \in
\Omega^0(\Gg_1 \timests \Gg_1)$, a 1-form $A \in
\Omega^1(\Gg_1)$, a 2-form $F \in \Omega^2(\Gg_0)$ and a 3-form
$K \in \Omega^3(\Gg_0)$ satisfying:
\begin{itemize}
\item $K=dF$
\item $\pi_2^*F - \pi_1^*F = dA$ and
\item $\pi_1^*A + \pi_2^*A +m^*i^*A = -\sqrt{-1} \theta^{-1} d\theta$
\end{itemize}
 A connection is {\it flat} if the curvature
$K$ vanishes.

\section{The Loop Groupoid}
 \label{sectionloopgroupoid}

\subsection{}We would like to define the loop groupoid associated to
a groupoid in a similar manner to the definition of the free loop
space associated to a manifold. Therefore we want the loop
groupoid to be a category whose objects are
${\mathrm{Hom}}(S^1,G)$ in the category of orbifolds. Let us
remember that the category of orbifolds is the same as the
category of proper \'{e}tale topological groupoids modulo Morita
equivalence. Therefore we must first assign a groupoid to $S^1$
and then consider all Morita equivalent groupoids. This amounts
in a concrete language to consider finer and finer covers of the
circle. Moreover, this allows us to give the loop groupoid a
natural topology.
We will do so in a very concrete manner by taking certain ``well
behaved'' open covers of the circle.

\begin{definition}
Let $\FF$ be the family of all finite sets of the unit circle
$S^1$. We identify $\FF$ with the family of all finite sets of
$(0,1]$ using the exponential map. In other words the set
$\{q_1,q_2, \dots, q_n,q_0\} \subset (0,1]$ (always written in
increasing order $q_1 < q_2 < \cdots < q_n < q_0$) will be
indentified with the ordered $n$-tuple of points in the circle
given by $\{ e^{2\pi q_k i} | k=1,\ldots,n,0\}$. To every element
$\{q_1,q_2, \dots, q_n,q_0\} \in \FF$ and every $\epsilon>0$ we
will associate a unique open cover of the circle
$\VV(\{q_k\},\epsilon)$ defined by the sets $V_i:=(q_i -
\epsilon,q_{i+1} + \epsilon)$. We denote by $\FF(\epsilon)$ the
family of all such covers. Notice that $1$ is always in $ V_0$.
Since  for every $\epsilon>0$, $\FF \cong \FF(\epsilon)$ as sets,
then $\FF(\epsilon)$ acquires a partial ordering (that we call
refinement) induced by that on $\FF$ given by inclusion of sets.
In addition we say that a cover $\VV(p_1,p_2, \dots,
p_m,p_0;\epsilon)$ is a refinement of $\VV(q_1,q_2, \dots,
q_n,q_0;\epsilon')$ when $ \epsilon \leq \epsilon' $ and
$\{q_1,q_2, \dots, q_n,q_0\} \subset \{p_1,p_2, \dots, p_m,p_0\}$
as sets. Given any two covers of types $\FF(\epsilon)$ and
$\FF(\epsilon')$ that are related (one being refinement of the
other) then there is a canonical embedding of the open sets of
the coarser cover into the ones of the finer. Finally we will
denote by $\WW$ the set of all covers $\VV(q_1,q_2, \dots,
q_n,q_0;\epsilon) \in \FF(\epsilon)$ for some $\epsilon$ so that
$\epsilon$ is small enough to make all triple intersections of
elements of $\VV(q_1,q_2, \dots, q_n,q_0;\epsilon')$, $V_i \cap
V_j \cap V_j = \emptyset$ . The set $\WW$ is called the set of
admissible covers of $S^1$. The remarks of this paragraph amount
to a definition of a partial order in $\WW$ called refinement.
This partial order has the property that given any two elements
$\VV(p_1,p_2, \dots, p_m,p_0;\epsilon)\in\WW$ and $\VV(q_1,q_2,
\dots, q_n,q_0;\epsilon')\in\WW$ then there is an element of
$\WW$ refining both at the same time (take for example
$\VV(\{q_1,q_2, \dots, q_n,q_0\}\cup\{p_1,p_2, \dots,
p_m,p_0\};\min(\epsilon,\epsilon'))\in\WW$.)
\end{definition}

 Now we will define a unique groupoid associated to the circle $S^1$ and a cover of
 it (and therefore one such groupoid for every element of $\WW$.)
Let $\{V_\alpha\}_{0 \leq \alpha \leq n}$ be such a cover of
$S^1$ and let $W=\{W_\beta\}_\beta$ be the pullback of this cover
under the map
$$\real \stackrel{e^{2 \pi i}}{\longrightarrow} S^1$$
and let's call by $\underline{\real^W}$ the groupoid associated
to this cover, in other words
$$\underline{\real^W}_1 := \bigsqcup_{(\beta_1,\beta_2)} W_{\beta_1} \cap W_{\beta_2}
\ \ \ \ \ \ \ \ \underline{\real^W}_0:=\bigsqcup_{\beta}
W_{\beta}.$$ If we call by $\{W_\alpha^n\}_n$ $n\in \integer$ the
inverse image of the set $V_\alpha$ by the exponential map, we
can define an action of $\integer$  in the following way:
\begin{eqnarray*}
 \underline{\real^W}_1 \times \integer&  \to & \underline{\real^W}_1 \\
 ((x,W_{\beta \gamma}),l) & \mapsto & (x+l,W_{\beta \gamma}^l)
\end{eqnarray*}
where $x \in W_{\beta \gamma}$ and $x+l \in W_{\beta \gamma}^l$.
\begin{definition}
Let $\Ss^1_W$ be the groupoid
$$\begin{array}{c}
\underline{\real^W}_1 \times \integer\\
\twodownarrows \\
\underline{\real^W}_0
\end{array}$$
with maps
$$s((x,W_{\beta \gamma}),l)=(x,W_{\beta}) \ \ \ \ \ \ \
t((x,W_{\beta \gamma}),l)=(x+l,W_{\gamma}^l)$$
$$e(x,W_\beta) = ((x,W_{\beta \beta}),0) \ \ \ \ \ \ \
i((x,W_{\beta \gamma}),l)=((x+l,W_{\gamma \beta}^l),-l)$$
$$m\left[((x,W_{\beta \gamma}),l),((x+l,W_{\gamma \sigma}^l),k) \right]=
((x,W_{\beta \sigma}),l+k)$$
\end{definition}
\subsection{}It is an easy exercise to check that $\Ss^1_W$ is Morita
equivalent to $S^1$ where its groupoid structure is the trivial
one, i.e. $S^1 \twoarrows S^1$ with the source and the target
maps equal to the identity. Moreover, for another  cover $W'$
defined in the same way as $W$ that refines it, there is a
natural morphism of groupoids $\rho^{W'}_{W}:\Ss^1_{W'} \to
\Ss^1_W$ which also turns out to be Morita.

\begin{definition}
For $\Gg$ a topological  groupoid and an open cover $W$ of the
circle, the loop groupoid $\Loop\Gg(W)$ associated to $\Gg$ and
the open cover $W$ will be defined by the following data:
\begin{itemize}
\item Objects ($\Loop\Gg(W)_0$):
Morphisms $\Ss_W^1 \to \Gg$
\item Morphisms ($\Loop\Gg(W)_1$): For two elements in $\Loop\Gg(W)_0$, say
$\Psi, \Phi : \Ss^1_W \to \Gg$ , a morphism (arrow) from $\Psi$
to $\Phi$ is a map $\Lambda : \underline{\real^W}_1 \times
\integer\to \Gg_1$ that makes the following diagram commute
      $$
       \xymatrix{
        \underline{\real^W}_1 \times \integer \ar[r]^\Lambda \ar[d]_{s\times t}
       & \Gg_1 \ar[d]^{s \times t}\\
       \underline{\real^W}_0 \times \underline{\real^W}_0 \ar[r]_{(\Psi_0 , \Phi_0)}
  & \Gg_0 \times \Gg_0 }
      $$
and such that for $r \in \underline{\real^W}_1 \times \integer$
\begin{eqnarray}  \label{eqnbetweenmorphisms}
\Lambda(r)=\Psi_1(r) \cdot \Lambda(et(r))= \Lambda(es(r)) \cdot
\Phi_1(r).
\end{eqnarray}
\end{itemize}
\end{definition}
The composition of morphisms is defined pointwise, in other
words, for $\Lambda$ and $\Omega$ with
       $$\xymatrix{
       \Psi \ar@/^/[r]^\Lambda  & \Phi \ar@/^/[r]^\Omega & \Gamma}
       $$
we set $$\Omega \circ \Lambda(es(r)) := \Lambda(es(r)) \cdot
\Omega(es(r))$$ and
 $$\Omega \circ \Lambda(r):= \Omega \circ \Lambda (es(r)) \cdot \Gamma(r)=\Psi(r) \cdot
\Omega \circ \Lambda (et(r))$$ The last equation comes from the
following set of equalities:
\begin{eqnarray*}
\Omega \circ \Lambda (es(r)) \cdot \Gamma(r) & = & \Lambda(es(r))
\cdot
\Omega(es(r)) \cdot \Gamma(r) \\
& = & \Lambda(es(r)) \cdot \Phi(r) \cdot \Omega(et(r)) \\
& = & \Psi(r) \cdot \Lambda(et(r)) \cdot \Omega(et(r))\\
& = & \Psi(r) \cdot \Omega \circ \Lambda(et(r))
\end{eqnarray*}

The groupoid $\Loop\Gg(W)$ can be given the compact-open topology
to make it a topological groupoid. To do this consider the
following identities satisfied by any arrow $\Lambda \colon \Psi
\to \Phi$
$$\Psi_1(r) = \Lambda(r) \cdot \Lambda(e(t(r)))^{-1},$$
$$\Psi_0(x) = s(\Lambda(e(x))),$$
$$\Phi_1(r) = \Lambda(e(s(r)))^{-1} \cdot \Lambda(r),$$
$$\Phi_0(x) = t(\Lambda(e(x))).$$
This identities imply that $\Lambda$ determines it source $\Psi$
and its target $\Phi$. Let $\Lambda^o = \Lambda \circ e$. Then
from the previous identities is easy to see that from the pair
$(\Psi_1,\Lambda^o)$ we can recover $\Lambda$. Therefore we can
see the set of all such $\Lambda$'s, namely $\Loop \Gg (W)_1$ as a
subspace of the space of continuous maps
$\mathrm{Map}(\underline{\real^W}_1,\Gg_1)\times\mathrm{Map}(\underline{\real^W}_0,\Gg_1)$.
In this way the spaces $\Loop \Gg(W)_1$ and $\Loop \Gg(W)_0$
inherit the compact-open topology, making the groupoid $\Loop
\Gg(W)$ into a topological one.

%
%

From now on we restrict our attention to admissible covers of the
circle. We would like that two morphisms $\Psi_i: \Ss^1_{W_i} \to
\Gg$ ($i=1,2$), be equivalent if there exists a refinement
$W\in\WW$ of $W_1\in\WW$ and $W_2\in\WW$ such that the following
diagram of morphisms commute

      $$\xymatrix{
      \Ss^1_{W} \ar[r]^{\rho^{W}_{W_1}} \ar[d]_{\rho^{W}_{W_2}} & \Ss^1_{W_1} \ar[d]^{\Psi_1}\\
      \Ss^1_{W_2} \ar[r]_{\Psi_2} & \Gg }
      $$
This can be  achieved using a limit construction. For $W'$ a
refinement of $W$ there is a natural monomorphism of topological
groupoids
$$\Loop\Gg(W) \hookrightarrow \Loop\Gg(W')$$
(here is where we use the partial order in $\WW$) this allows us
to define
\begin{definition}
The loop groupoid $\Loop\Gg$ of the topological groupoid $\Gg$ is
defined as the monotone union (colimit) of the groupoids
$\Loop\Gg(W)$ where $W$ runs over the set $\WW$ of admissible
covers
$$\Loop\Gg := \lim_{\overrightarrow{W\in\WW}} \Loop\Gg(W).$$
\end{definition}

In this way the loop groupoid is naturally endowed with a
topology, becoming a topological groupoid.

\begin{rem}
Unless stated otherwise the groupoids we will focus our attention
on will be orbifolds (smooth, \'etale and proper groupoids)
although some of the results that follow can be generalized.
\end{rem}

Now let's see a property of the arrows in the loop groupoid of an
\'{e}tale groupoid $\Gg$. This paragraph justifies us to
interpret the loop groupoid $\Loop \Gg$ of an (orbifold)
\'{e}tale proper groupoid $\Gg$ as an infinite dimensional
orbifold groupoid.

\begin{lemma} \label{arrowdeterminedbypoint}
Let $\Gg$ be an \'{e}tale groupoid. An arrow $\Lambda\in \Loop
\Gg_1$ between $\Psi$ and $\Phi$ is determined completely by
$\Psi$ and the image $\Lambda^o(x)$  of a marked point $x \in
\underline{\real^W}_0$. Here $\Lambda^o = \Lambda \circ e $ is
the restriction of $\Lambda$ to the identity arrows of $\Gg$.
\end{lemma}
\begin{proof}
First is clear that $\Lambda$ is completely determined by
$\Lambda^o$ because $$ \Lambda(r) = \Psi(r) \cdot
\Lambda^o(t(r)).$$ What this lemma is saying is that the image of
any other point $y \in \underline{\real^W}_0$ is uniquely
determined by $\Lambda^o(x)$ and the morphism $\Psi$. Let's
suppose without loss of generality that $x \in V_1$ and $y \in
V_2$, both sets in $W$ such that $V_{12}:=V_1 \cap V_2 \neq
\phi$. As $\Gg$ is \'{e}tale the image of $V_1$ under $\Lambda$
is determined by the point $\Lambda^o(x)$ and $\Psi_0(V_1)$; this
because we can find a neighborhood $U_1$ of $\Lambda^o(V_1)$ in
$\Gg_1$ homeomorphic to $s(U_1)$ which contains $\Psi_0(V_1)$.
Now if $z \in U_{12}$, labeling it $z_1$ when we see it in $V_1$
and $z_2$ when in $V_2$, then
$$\Lambda^o(z_2) = \Psi_1(z)^{-1} \Lambda^o(z_1) \Psi_1(z).$$
Applying the same procedure as before, the set $\Lambda(V_2)$ is
determined by $\Psi_0(V_2)$ and $\Lambda^o(z_2)$; and hence it
determines $\Lambda^o(y)$.
\end{proof}

This allows us to get the following:

\begin{cor}
Let $\Lambda_0$ and $\Lambda_1$ be arrows joining $\Psi$ and
$\Phi$. If $\Lambda_0(x)=\Lambda_1(x)$ then $\Lambda_0=\Lambda_1$.
\end{cor}

\begin{cor}
If $\Gg$ is an \'{e}tale proper groupoid, then the loop groupoid
$\Loop \Gg$ is \'{e}tale and the local isotropy at any object
$\Psi$ is finite.
\end{cor}
\begin{proof}
The arrows between two morphisms are determined by the image of
only one point $\Lambda(x)$ in $\Gg_1$. If the source and the
target of $\Lambda$ are equal, then the same is true for the
arrow $\Lambda(x)$. As the isotropy in $\Gg$ is finite, the
result follows. This implies that the source map of the loop
groupoid is \'{e}tale.
\end{proof}

We say then that $\Loop \Gg$ represents an infinite dimensional
orbifold.

\subsection{} \label{morphismtrivialcover}
In order to motivate this definitions  let's consider the case
when $\underline{\real}_1=\underline{\real}_0=\real$. Let
$\Psi:\Ss^1 \to \Gg$ be a morphism, then it is determined by the
image of $\real \times \{1\}$ under $\Psi_1$; clearly
$\Psi_1(x,0)=es(\Psi_1(x,1))$, $\Psi_0(x)=s(\Psi_1(x,1))$ and
$\Psi_0(x+1)= t(\Psi_1(x,1))$, also using the fact that
$$(x,n)=(x,1)\cdot(x+1,1)\cdot \dots \cdot (x+n-1,1) $$
we obtain
$$\Psi_1(x,n) =\Psi_1(x,1)\cdot \Psi_1(x+1,1)\cdot \dots \cdot \Psi_1(x-n+1,1).$$
The use of the word loop is justified because the function $\Psi_0
| _{[x,x+1]}$ gives us a path from  the source to the target of
$\Psi_1(x,1)$
      $$
       
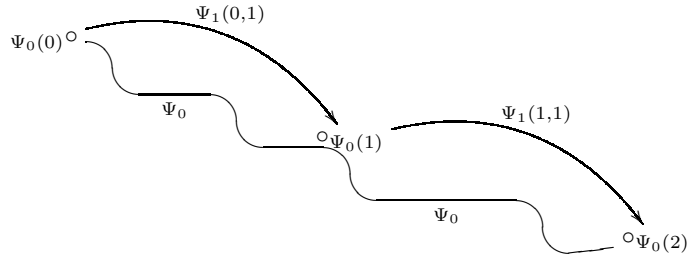
\xymatrix{
       {}_{\Psi_0(0)}{\circ} \ar@{-} `r[d] `[rr] `r[drr]_{\Psi_0} `[drrr]
       `r[ddrrr] `[drrrrr] `r[ddrrrrr]_{\Psi_0} `[ddrrrrrr] [ddrrrrrr] \ar@/^2pc/[drrr]^{\Psi_1(0,1)}
       & & & & & &\\
       & &  &{\circ}_{\Psi_0(1)} \ar@/^2pc/[drrr]^{\Psi_1(1,1)} & & &\\
       & & & & & & \circ_{\Psi_0(2)} }
      $$
and when we consider the coarse moduli space of the groupoid the
path determined by $\Psi_0 | _{[0,1]}$ becomes a loop. The
morphisms of $\Loop\Gg$ relate the loops that in the coarse
moduli space are identified, so for $\Lambda$ an arrow from
$\Psi$ to $\Phi$ we
 have the following diagram
      $$
       \xymatrix{
       {}^{\Psi_0(x)}{\circ} \ar@{-} `r[d] `[rr] `r[drr]_{\Psi_0} `[drrr] [drrr]
        \ar@/^1pc/[drrr]^{\Psi_1(x,1)} \ar@/^2pc/[rrrr]^{\Lambda(x,0)}
       & & & &{\circ}^{\Phi_0(x)}\ar@{-} `r[d] `[rr] `r[drr]_{\Phi_0} `[drrr] [drrr]
       \ar@/^1pc/[drrr]^{\Phi_1(x,1)} & & &\\
       & &  &{\circ}_{\Psi_0(x)} \ar@/_1pc/[rrrr]^{\Lambda(0,x+1)} & &  & & {\circ}_{\Phi_0(x+1)}}
      $$
\subsection{} Now we set to the task of showing that $\Loop \Gg$
doesn't depend on the Morita equivalence class of $\Gg$. In
particular $\Ll$ is a functor from groupoids to groupoids.
\begin{definition}
A morphisms of groupoids $\Ff: \Hh \to \Gg$ induces naturally a
morphism between the loop groupoids $\Loop \Ff : \Loop \Hh \to
\Loop \Gg$  in the following way:
\begin{itemize}
\item Objects:
      $$
       \xymatrix{
       (\Loop \Hh)_0   \ar[r]^{(\Loop \Ff)_0} & (\Loop \Gg)_0 }$$
       $$\xymatrix{
       \ \ \Psi \ \ \ar@{|->}[r] & \Ff \circ \Psi}
      $$
\item Morphisms:  For $\Lambda \in \Loop \Hh$, an arrow between $\Psi$ and $\Phi$,
we define  $(\Loop \Ff)_1 (\Lambda) := \Ff_1 \circ \Lambda$ as an
arrow between $(\Loop \Ff)_0(\Psi)$ and $(\Loop \Ff)_0(\Phi)$
     $$
      \xymatrix{
       \integer \times \underline{\real}_1 \ar[r]_\Lambda \ar[d]_{s\times t}
      \ar@/^1pc/[rr]^{\Ff_1\circ\Lambda}
      & \Hh_1 \ar[d]^{s \times t} \ar[r]_{\Ff_1} & \Gg_1 \ar[d]^{s \times t} \\
      \underline{\real}_0 \times \underline{\real}_0 \ar[r]^{(\Psi_0 , \Phi_0)}
       & \Hh_0 \times \Hh_0
      \ar[r]^{(\Ff_0,\Ff_0)} & \Gg_0 \times \Gg_0 }
     $$
\end{itemize}
The morphism $\Loop \Ff$ is continuous because the  limit
commutes with the composition of $\Ff$.
\end{definition}
\begin{lemma} \label{moritaloopgroupoids}
If $\Ff : \Hh \to \Gg$ is a Morita morphism of \'{e}tale
groupoids then the induced morphism $\Loop \Ff : \Loop \Hh \to
\Loop \Gg$ is also Morita.
\end{lemma}
\begin{proof}
First we need to check that the following square is a fibered
product
       $$
        \xymatrix{
        (\Loop \Hh)_1  \ar[d]^{s \times t} \ar[rr]^{(\Loop \Ff)_1} & & (\Loop \Gg)_1
        \ar[d]^{s \times t} \\
        (\Loop \Hh)_0 \times (\Loop \Hh)_0 \ar[rr]^{((\Loop \Ff)_0,(\Loop \Ff)_0)} & &
        (\Loop \Gg)_0 \times (\Loop \Gg)_0 }
       $$
So for $\Psi, \Phi \in (\Loop \Hh)_0$ and $\Omega \in (\Loop
\Gg)_1$ with
$$s(\Omega)= \Ff \circ \Psi \ \ \ \ \ \ \ \  \ \ \ \ \ t(\Omega) = \Ff \circ \Phi$$
we need to show that there is only one $\Lambda \in (\Loop
\Hh)_1$ such that $(\Loop \Ff)_1(\Lambda) = \Omega$,
$s(\Lambda)=\Psi$ and $t(\Lambda)=\Phi$.
 Let's work pointwise; for $a \in e(\underline{\real}_0)$
let
$$r:=\Omega(a) \ \ \ \ \ \ x:= s(\Psi_1(r)) \ \ \ \ \ \ y:= s(\Phi_1(r))$$
then $s(r) = \Ff_0(x)$ and $t(r) = \Ff_0(y)$; as $\Ff$ is a
Morita morphism then there is only one element $w \in \Hh_1$ such
that $s(w)=x$, $t(w)=y$ and $\Ff_1(w)=r$; then $\Lambda(a):=w$.
So we have a well defined map $\Lambda: \underline{\real}_1
\times \integer \to \Hh_1$ with the required properties. As the
maps $\Ff_0$ and $\Ff_1$ are \'{e}tale then is easy to see that
$\Lambda$ is continuous. We still have to prove that for $\alpha
\in \underline{\real}_1 \times \integer$ we have
$$\Lambda(es(\alpha)) \cdot \Phi_1(\alpha) = \Psi_1(\alpha) \cdot \Lambda(et(\alpha))$$
Let $v:=\Psi_1(\alpha) \cdot \Lambda(et(\alpha))\cdot
\Phi(\alpha)^{-1} \cdot \Lambda(es(\alpha))^{-1}$, as $\Omega$
belongs to $(\Loop \Gg)_1$ then $\Ff_1(v)=es(\Ff_1(v))$, but
$es(v)$ also maps to the same morphism under $\Ff_1$; then by the
properties of the Morita morphism we get that $v = es(v)$.

Now we'll see that the map $\Loop \Ff_0 : \Loop \Hh_0 \to \Loop
\Gg_0$ is surjective. Let $W$ be a cover associated to the tuple
$\{q_1, \dots, q_n, q_0\}$ and $\Psi : \Ss^1_W \to \Gg$ an object
in $\Loop\Gg$. Let $V_i=(q_i - \epsilon, q_{i+1} +\epsilon)$ be
an  open set of $W$ and $K_i \subset \Gg_0$ the image under
$\Psi_0$ of the compact set $[q_i,q_{i+1}]$. As the map $\Ff_0$
is \'etale, there are a finite number of open sets $U^i_j \subset
\Hh_0$, $0 \leq j \leq l_i$  diffeomorphic to their image under
$\Ff_0$, i.e. $U^i_j \cong \Ff_0(U^i_j)$, such that the sets
$\Ff_0=(U^i_j)$ cover $K_i$ ($K_i$ is compact). This is done in
such a way that we can partition the interval $[q_i,q_{i+1}]$ in
$l_i+1$ pieces
  $q_i=p^i_0< p^i_1< \cdots<p^i_{l_i}=q_{i+1}$
such that
$$\Psi_0([p^i_j,p^i_{j+1}]) \subset \Ff_0(U^i_j).$$
We do the same for $i$ in $0 \leq i \leq n$ and  we take the
cover $W'$ associated to the tuple given by the $p^i_j$'s. Define
$\Phi_0 : \underline{\real^{W'}}_0 \to \Hh_0$ by
\begin{eqnarray*}
 \Phi_0 : (p^i_j -\epsilon' ,p^i_{j+1}+\epsilon') & \to & \Hh_0 \\
x & \mapsto & (\Ff_0|_{\Ff_0(U^i_j)})^{-1} \circ \Psi_0 \circ
\rho^{W'}_W(x).
\end{eqnarray*}
The construction of $\Phi_1$ is straight forward using the fact
that the map $\Ff$ is a Morita equivalence, i.e. the image of
$\Phi_0$ and the cartesian square of \ref{Moritasquare} defines
it uniquely. Even though $\Ff \circ \Phi$ is not equal to $\Psi$,
it is clear  that both are in the same class of  the limit.

\end{proof}

\subsection{} We can see that the loop groupoid defined in this
way is an immense object. We want this definition to be more
manageable, so in some cases, depending on whether  our groupoid
$\Gg$ has some additional structure we will see, in the following
sections, that we can considerably simplify the construction.
Notice as well that it is closely related to the path
constructions of \cite{BridsonHaefliger, Mrcun}

\subsection{} The loop groupoid $\Loop \Gg$ has a natural action of
$\real$ (that actually descends to an action of $S^1$ given by
rotating the loop). Let's recall that the groupoid $\Ss^1$ is
subordinated to a cover $W$ of $\real$; for $z \in \real$, let
$W^z$ be the cover of $\real$ obtained from $W$ after shifting by
$z$. To be more explicit $W^z =\{V^z | V \in W\}$ where
$V^z=\{x+z|x \in V\}$. Using this new cover we can define the
action of $\real$ into $\Loop \Gg$.

\begin{definition} For $z \in \real$ and $\Psi \in (\Loop \Gg)_0$, let
$\Psi^z$ be the morphism
$$\Psi^z : \Ss^1_{W^z} \to \Gg $$
defined in the following way
$$\Psi^z_1 ((x,V_{\alpha \beta}^z),l) := \Psi ((x-z,V_{\alpha \beta}),l)$$
$$\Psi^z_0(x,V_\alpha^z):=\Psi_0(x-z,V_{\alpha})$$
Obtaining an action of $\real$ into $\Loop \Gg$
\begin{eqnarray*}
(\Loop \Gg)_0 \times \real & \stackrel{z}{\to} & (\Loop \Gg)_0\\
(\Psi,z) & \mapsto & \Psi^z
\end{eqnarray*}
\end{definition}

We will denote by $\Loop \Gg ^\real$ the subgroupoid of $\Loop
\Gg$ that is fixed under the action of $\real$. The great
generality that we are allowing in our definition of the loop
groupoid may persuade us that $\Loop \Gg^\real$ is very
complicated groupoid. This is not the case amd we identify it now.
Every morphism in $\Loop \Gg ^\real$ is constant when is seen in
the coarse moduli space of $\Gg$, but in $\Gg$ several objects map
into that constant point. The groupoid $\Loop \Gg^\real$ is
taking into account all that information. The following technical
lemmas will make more explicit the idea just mentioned.

\begin{lemma}
Let $\Psi$ be an object in $\Loop \Gg^\real$, then $\Psi_1$ and
$\Psi_0$ are locally constant.
\end{lemma}
\begin{proof}
Let $V \in W$ and $z \in \real$ such that $V \cap V^z \not= \phi$.
As $V\cap V^z$ belongs to a common refinement of $W$ and $W^z$ and
$\Psi=\Psi^z$ then we must have that
$$\Psi_0(x,V) = \Psi^z_0(x,V^z) = \Psi_0(x-z,V)$$
where the first equality is the condition of being a fixed element
under the action of $\real$, and the second is just the
definition of $\Psi^z$. Then we can see that $\Psi_0$ is constant
in all $V$. For $\Psi_1$ the same argument is applied.
\end{proof}

\begin{lemma}
 \label{lemmacambiocubiertas}
Let $\Psi$ be in $\Loop \Gg^\real$, then there exist another
object $\Phi$ in $\Loop \Gg^\real$  with
$$\Phi : \Ss^1_\real \to \Gg$$
 defined over the trivial cover of $\real$ such that there is
an arrow $\Lambda$ joining $\Psi$ and $\Phi$.
\end{lemma}
\begin{proof}
The idea underlying in this lemma is the fact that all the
objects in
 the image of $\Psi_0$ are related by arrows in $\Gg_1$.\\
Let's fix the point $(0,V_0)$ in $\underline{\real}_0$ (as we
said before, $V_0$ always conatins the marked point of the
circle); for $(y,U)$ another point in  $\underline{\real}_0$ there
exist open sets $V=U_0, U_1, \dots, U_n=U$ in the cover such that
$U_i \cap U_{i+1} \not= \phi  $; let $y_i$ be any point in
$U_{i,i+1}:=U_i \cap U_{i+1}$. As $\Psi$ is a locally constant
morphism $$\Psi_0(y_i,U_i)=\Psi_0(y_{i+1},U_i)$$ and
$\Psi_1((y_i,U_{i,i+1}),0)$ is an arrow between
$$ \Psi_0(y_i,U_i) \ \ \ \ \ \mbox{    and   } \ \ \ \ \ \ \Psi_0(y_i, U_{i+1})$$
Thus
\begin{eqnarray} \label{Ar}
Ar^{x,V}_{y,U} := \Psi_1((y_0,U_{0,1}),0) \cdots
\Psi_1((y_{n-1},U_{n-1,n}),0)
\end{eqnarray}
is an arrow in $\Gg_1$ joining $\Psi_0(x,V)$ and $ \Psi_0(y,U)$.
It is worth pointing out that because of the properties of the
groupoids the map $Ar^{x,V}_{y,U}$ is independent of the chosen
sets $U_i$; is also clear that $Ar$ only depends on the open sets
$V$ and $U$, so we can suppress the points $x$ and $y$;
$Ar^V_U=Ar^{x,V}_{y,U}$ . Let $\Phi: \Ss^1_W \to \Gg$ be the
morphism defined in the following way (as the morphism is locally
constant we will define it on the open sets):
$$\Phi_0(U) :=\Psi_0(V)$$
\begin{eqnarray} \label{equationofshifting}
\Phi_1(U_{\alpha\beta},l) := Ar^V_{U_\alpha} \cdot
\Psi_1(U_{\alpha\beta},l) \cdot \big(Ar^V_{U^l_\beta} \big)^{-1}
\end{eqnarray}
Replacing \ref{Ar} when $l=0$ we get that
$$\Phi_1(U_{\alpha\beta},0) = \Psi_1(V,0) = e\Psi_0(V)$$
and from these it follows that
$$\Phi_1((U_{\alpha\beta},l)=\Psi_1(V,l).$$
So we can make $\Phi$ to be defined on the trivial cover of
$\real$, $\Phi:\Ss^1_\real \to \Gg$
$$\Phi_0(y) :=\Psi_0(V) \ \ \ \ \ \mbox{and} \ \ \ \ \  \ \Phi_1(y,l):= \Psi_1(V,l)$$

Now we need to construct an arrow $\Lambda$ between $\Psi$ and
$\Phi$, and is very clear how it should be defined:
$$\Lambda(U_{\alpha\alpha},0) :=Ar^V_{U_\alpha}$$
then equation \ref{equationofshifting} gives explicitly the
condition that $\Lambda$ has to fulfill.
\end{proof}
So, for the groupoid $\Loop \Gg^\real$ we can restrict our
attention to the morphisms $\Phi : \Ss^1_\real \to \Gg$ defined
on the trivial cover of $\real$,
 but we just proved
that they are locally constant; so we get
\begin{theorem}
The groupoid $\Loop \Gg^\real$ is Morita equivalent to the
groupoid $Hom(\overline{\integer}, \Gg)$ (with the natural
topology) whose objects are morphisms $\phi: \overline{\integer}
\to \Gg$ with $s\phi(1)=t\phi(1) = \phi(0)$ and whose arrows
$\phi \stackrel{\lambda}{\to} \psi$ are maps $\lambda : \integer
\to  \Gg_1$ with
$$\lambda(0) \cdot \psi_1(n) =\phi_1(n) \cdot \lambda(0)=\lambda(n)$$
\end{theorem}
\begin{proof}
From $\Phi: \Ss^1_\real \to \Gg$ obtained in the previous lemma,
which is locally constant, we can define $\phi:
\overline{\integer} \to \Gg$ as $\phi_1(n):=\Phi_1(x,n)$ and
$\phi_0:= \Phi(x)$ for any $x$ in $\real$. From the fact that
$s\Phi_1(0,1)= \Phi_0(0)=\Phi_0(1)=t\Phi_1(0,1)$ we get that
$s\phi(1)=t\phi(1) = \phi(0)$.

For an arrow $\Phi \stackrel{\Lambda}{\to} \Psi$, with $\Psi$
another object in $\Loop \Gg^\real$ as in previous lemma, we can
define $\lambda(n):=\Lambda(x,n)$ and clearly we get what we
needed.

The Morita equivalence follows from the fact that we are getting
rid
 of extra objects and arrows, the proof is straight forward.
 \end{proof}
This description of  $\Loop \Gg^\real$ matches the one of what is
known as the {\it inertia groupoid}.
\begin{definition}
The inertia groupoid $\wedge \Gg$ is defined in the following way:
\begin{itemize}
\item Objects $(\wedge \Gg)_0$: Elements $v \in \Gg_1$ such that $s(v) = t(v)$.
\item Morphisms $(\wedge \Gg)_1$: For  $v,w \in (\wedge \Gg)_0$ an
arrow $v \stackrel{\alpha}{\to} w$ is an element $\alpha \in
\Gg_1$ such that $v \cdot \alpha = \alpha \cdot w$
        $$
         \xymatrix{
         \circ \ar@(ul,dl)[]|{v} \ar@/^/[rr]|{\alpha}
         &&\circ \ar@(dr,ur)[]|{w^{-1}} \ar@/^/[ll]|{\alpha^{-1}}
         }$$
\end{itemize}
It is known that the inertia groupoid in the case of an orbifold
matches with what is commonly known in the literature by twisted
sectors (see \cite{Ruan}),
 thus this is
a natural way to define them. We conclude this section with the
following statements:

\end{definition}

\begin{proposition} The groupoid $\Loop \Gg^\real$ is Morita
equivalent to the inertia groupoid $\wedge \Gg$.
\end{proposition}
\begin{proof}
Follows from the fact that the groupoids $\wedge \Gg$ and
$Hom(\overline{\integer}, \Gg)$ are clearly diffeomorphic
(homeomorphic).
\end{proof}
And by using the result of the appendix we get.

\begin{cor}
For an orbifold $\Gg$, the groupoid $\Loop \Gg^\real$ is Morita
equivalent to the twisted sectors $\widetilde{\Sigma_1 \Gg}$ of
$\Gg$.
\end{cor}

\section{Finite Groups and Lie Groups.}

\subsection{}Even for the case in which our groupoid $\Gg$ happens to
be a finite group the theory is quite interesting. For $G$ a
finite group acting on the connected space $X$ and $\Gg$ the
groupoid associated to the global quotient $[X/G]$ one does not
need to consider the morphisms from $\Ss^1$ to $\Gg$ with all
possible open covers of $\real$. As the space $X$ is connected,
for any $\Psi : \Ss^1_W \to \Gg$ there exist another morphism
$\Phi : \Ss^1_\real \to \Gg$ subordinated to the trivial cover of
the reals (i.e. $\underline{\real}_1=\underline{\real}_0=\real$)
and an arrow in $(\Loop \Gg)_1$ relating the two.

 The idea here
is that we are getting rid of unnecessary arrows as well as
objects without changing the relevant information of the
groupoid, so we have
\begin{lemma}  \label{lemmaglobalquotient}
Let $\Gg=[X/G]$ be the groupoid associated to the global quotient
$X/G$ for $G$ a finite group and $X$ a  $G$ space, then for the
loop groupoid it suffices to take morphisms from $\Ss^1$ to $\Gg$
associated to the trivial cover of $\real$.
\end{lemma}
\begin{proof}
The argument follows the lines of the proof of lemma
\ref{lemmacambiocubiertas} although in this case is a bit more
complicated. For every $\Psi \in \Loop [X/G]$ we are going to
associate a morphism subbordinated to the trivial cover as
follows. Let $W=\{V_i\}_{0 \leq i \leq n}$ be the cover of the
circle on which $\Psi$ is defined and   $\{U_i^j\}_{0 \leq i \leq
n, j \in \integer}$ the cover of $\real$ given by the pullback of
$W$ under the exponential map (in the same way that is done at
the beginning of chapter \ref{sectionloopgroupoid}). To make the
notation less cumbersome we will write $U_{i,i+1}^j:=U_i^j \cap
U_{i+1}^j$ and $U_{n,0}^{j,j+1}:=U_n^j \cap U_0^{j+1}$. We will
define another morphism $\Psi'$ subbordinated to the same cover
$W$ that will differ to $\Psi$ only in the information relevant
to the open sets $\{U_1^j\}_{j \in \integer}$ (we want that
$\Psi'_0$ could be defined on $U_0^j \cup U_1^j$). We are going
to abuse notation and we will define $\Psi'$ not point by point
but by defining the images of the sets $U_i^j$'s.

Let $g_j \in G$ be such that $\Psi_1(U_{0,1}^j,0) \subset X
\times \{g_j\}$, and $\pi_i :X\times G \to X, G$ the projections
on the first and second coordinates respectively. Let's define
$\Psi'$ as follows:
\begin{eqnarray*}
\Psi'_1(U_i^j,k) & :=& \Psi_1(U_i^j,k) \mbox{ for } i \neq 1, j,k \in \integer\\
\Psi'_0(U_1^j) &:=& \Psi_0(U_1^j)\cdot g_j^{-1} \\
\Psi'_1(U_{0,1}^j,0) &:=& (\pi_1(\Psi'_1(U_{0,1}^j,0)), id) \\
\Psi'_1(U_{1,2}^j,0) &:=& (\pi_1(\Psi'_1(U_{1,2}^j,0)) \cdot g_j^{-1}, g_j \pi_2(\Psi'_1(U_{1,2}^j,0)) )\\
\Psi'_1(U_1^j,1) &:=& (\pi_1(\Psi'_1(U_1^j,1)) \cdot g_j^{-1},
g_j \pi_2(\Psi'_1(U_{1,2}^j,0)) g_{j+1}^{-1}).
\end{eqnarray*}
The rest of the maps are determined by the previous ones. Now, as
in $\Psi'$ the image of $\Psi'_1(U_{0,1}^j,0)$ lies in $X \times
\{id\}$ we could say that $\Psi'$ is subbordinated to the cover
 $W' = W - \{V_0,V_1\} \cup \{V_0 \cup V_1\}$.

We could do this process $n$ times and we obtain a morphism
$\Psi'$ that is subbordinated to the cover given by the sets $U^j
: = U_0^j \cup U_1^j \cup \cdots \cup U_n^j$ (note that $U^j = (j
-\epsilon, j+1 +\epsilon)$). As before let's denote the set $U^j
\cap U^{j+1}$ by $U^{j,j+1}$. Let $k_i,h_i \in G$ be such that
$\Psi'_1(U^{i-1,i},0) \subset X \times k_i$ for $i>0$ and
$\Psi'_1(U^{-(i-1),-i},0) \subset X \times h_{i}$ for $i>0$; we
can define $\bar{\Psi}$ as follows:
\begin{eqnarray*}
\bar{\Psi}_0(U^0) & :=& \Psi'_0(U^0)\\
\bar{\Psi}_0(U^i) & :=& \Psi'_0(U^i) k_i^{-1}  k_{i-1}^{-1} \cdots k_1^{-1} \mbox{ for } i >0\\
\bar{\Psi}_0(U^i) & :=& \Psi'_0(U^i) h_{-i}^{-1}  h_{-(i+1)}^{-1} \cdots h_1^{-1} \mbox{ for } i <0\\
\bar{\Psi}_1(U^{i-1,i},0) & :=& (\pi_1(\Psi'_1(U^{i-1,i},0)) k_i^{-1}  k_{i-1}^{-1} \cdots k_1^{-1}, id) \mbox{ for } i >0\\
\bar{\Psi}_1(U^{i+1,i},0) & :=& (\pi_1(\Psi'_1(U^{i+1,i}0))
h_{-i}^{-1}  h_{-(i+1)}^{-1} \cdots h_1^{-1}, id)
 \mbox{ for } i <0\\
\bar{\Psi}_1(U^{i},1) & :=& (\bar{\Psi}_0(U^i),
k_1 \cdots k_i  \pi_2(\Psi'_1(U^{i},1))   k_i^{-1} \cdots  k_1^{-1}  ) \mbox{ for } i \geq0\\
\bar{\Psi}_1(U^{i},-1) & :=& (\bar{\Psi}_0(U^i), h_1 \cdots h_{-i}
\pi_2(\Psi'_1(U^{i},-1))  h_{-i}^{-1} \cdots h_1^{-1}) \mbox{ for } i \leq0\\
\end{eqnarray*}
In this way $\bar{\Psi}$ is subbordinated to the trivial cover of
the circle and is clear by lemma \ref{arrowdeterminedbypoint}
that there is only one arrow $\Lambda$ between $\Psi$ and
$\bar{\Psi}$ (note that for $0 \in U^0_0$,
$\Psi_0(0)=\bar{\Psi}_0(0)$).

\end{proof}

\subsection{} For example if $\Gg=\overline{G}$ for $G$ a finite groupoid
(i.e. $X = *$) and a morphism $\Phi_1 : \real \times \integer \to
G$, $\Phi_0: \real \to *$ the maps $\Phi_0$ and $\Phi_1(n,\_)$
are all constant. So we only need to consider the morphisms of
the groupoid $\overline{\integer}$ to $\overline{G}$, which is
precisely $Hom(\integer,G)$. It's also easy to see that for any
$h \in G$ there is an arrow between $\rho \in Hom(\integer,G)$
and $\rho^h \in Hom(\integer,G)$ where $\rho^h(n) = h^{-1}
\rho(n) h$. Hence

\begin{lemma}
 \label{moritafinitegroup}
Let $G$ be a finite group. Then the loop groupoid $\Loop
\overline{G}$ of $\overline{G}$ is Morita equivalent to the
groupoid $[Hom(\integer,G) / G]$ where $G$ acts on
$Hom(\integer,G)$ by conjugation.
\end{lemma}

So $\Loop \overline{G}$ can also be seen as the groupoid $G\times
G \twoarrows G$ where $s(g,h)=g$ and $t(g,h) = h^{-1}gh$, in
other words $G$ acts on itself by conjugation. It is clear that
the groupoids
\begin{equation}\label{moritaeq1}
\begin{array}{c}
G\times G\\
\twodownarrows \\
G
\end{array}
\stackrel{M}{\cong} \bigsqcup_{(g)} \left(
\begin{array}{c}
* \times C(g)\\
\twodownarrows\\
*
\end{array} \right)
\end{equation}
  are Morita equivalent, where the right hand side runs over the conjugacy classes
of elements in G and $C(g)$ is the centralizer of $g$.

\begin{rem}
The groupoid $\bigsqcup_{(g)} \overline{C(g)}$ is also a
subgroupoid of $\Loop \overline{G}=\wedge \overline{G}$.
$$
\bigsqcup_{(g)} \left(
\begin{array}{c}
* \times C(g)\\
\twodownarrows\\
*
\end{array} \right)
\hookrightarrow
\begin{array}{c}
G\times G\\
\twodownarrows \\
G
\end{array}
$$
\end{rem}

 Now applying the classifying space functor to the map of lemma \ref{moritafinitegroup} we obtain
\begin{lemma} For $G$ a finite group
$$B \Loop \overline{G} \simeq \bigsqcup_{(g)} BC(g).$$
\end{lemma}
And using  a result in algebraic topology \cite[p. 512]{Adem}
that says that there exists a natural map
$$f :  \bigsqcup_{(g)} BC(g)\simeq \LL B G$$ which is a homotopy equivalence,  get that
\begin{proposition}
For $G$ a finite group, there is a natural homotopy equivalence
$$B \Loop \overline{G} \simeq \LL BG.$$
\end{proposition}

\subsection{}\label{sectiongerbeongroup}
 Let's consider now a gerbe over the groupoid $\overline{G}$; from
\cite[Ex. 6.1.4]{LupercioUribe} we know it consists of a line
bundle $\LLL$ over $G$ and a 2-cocycle $\theta : G \times G \to
\U{1}$ such that
$$\LLL_g^{-1} \cong \LLL_{g^{-1}} \ \ \ \ \ \mbox{and} \ \ \ \ \ \
\LLL_g\LLL_h \stackrel{\theta(g,h)}{\cong} \LLL_{gh}.$$ Using
this information we would like to define a line bundle on the
inertia groupoid
 as follows. The arrows in the inertia groupoid of $\overline{G}$ relate
loops via conjugation. Let $g,h \in G$ where $g$ represents a map
$\overline{\integer} \to \overline{G}$ and $h$  an arrow in
$(\Loop \overline{G})$ joining $g$ and $h^{-1}gh$.

          $$
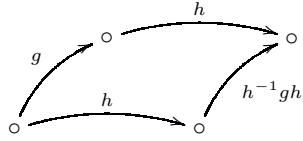
\xymatrix{
         & \circ \ar@/^/[rr]^{h}  & & \circ\\
         \circ \ar@/^/[ur]^g \ar@/^/[rr]^h & & \circ \ar@/^/[ur]_{h^{-1}gh}& }$$

we have
$$\LLL_g \LLL_h \stackrel{\frac{\theta(g,h)}{\theta(h,(h^{-1}gh))}}{\cong}
 \LLL_{h} \LLL_{(h^{-1}gh)} $$
then  the map
$$\rho(g,h) = \frac{\theta(g,h)}{ \theta(h,h^{-1}gh)}$$
 relates the fibers $\LLL_g$ and $\LLL_{h^{-1}gh}$. What is remarkable about this map $\rho$ is that
it gives a morphism from $\wedge \overline{G}$ to
$\overline{\U{1}}$.

\begin{lemma} \label{lemmalinebundlediscretetorsion}
Let $\wedge \overline{G} = G \times G \twoarrows G$ be the
inertia groupoid of $\overline{G}$, then
$$\rho :\wedge \overline{G} \to \overline{\U{1}}$$
is a morphism of groupoids. Moreover, it produces a line bundle
on $\wedge \overline{G}$ which matches the inner local system of
\cite{Ruan}.
\end{lemma}
\begin{proof} See lemma \ref{lemmamorphismofgroupoids}.
\end{proof}

Using the subgroupoid representation $ \sqcup_{(g)}
\overline{C(g)} \hookrightarrow
 \wedge \overline{G}$  we can consider
the map $\rho$ restricted to the centralizers $C(g)$, we get that
\begin{proposition}\label{rep1}
The map $$\rho(g,\cdot) : C(g) \to \U{1}$$ is a representation
that defines a line bundle on the groupoid $\overline{C(g)}$.
 This description matches precisely the inner local system defined by
Ruan \cite{Ruan} that comes from a discrete torsion
 $\theta \in H^2(G,\U{1})$.
\end{proposition}

So we can give now a natural definition of some concepts
introduced in \cite{Ruan},

\begin{proposition}
The inner local system associated to a discrete torsion $\theta
\in H^2(G,\U{1})$ arise from the relation on the fibers of the
gerbe associated to $\theta$ over constant loops in $\Loop
\overline{G}$, given by the arrows of $(\Loop
\overline{G})^\real$ as in \ref{rep1}.
\end{proposition}

\subsection{} Let us now briefly consider now the case in which $G$ is a
\emph{Lie group}. The groupoid $ \Hh = (G\times G
\rightrightarrows G)$ (with $G$ acting on $G$ by conjugation,) is
also naturally inside $\Loop\overline{G}$.

Consider the trivial $G$-principal bundle $P=G\times S^1 \to S^1$
over the circle. Let $\Aa$ be the space of connections on $P$,
and let $\GG$ be the gauge group of automorphisms of $P$.

Let $\Gg$ be the groupoid $\Aa \times \GG \rightrightarrows \Aa$,
where $s(A,g)=A$ and $t(A,g)= g^* A$.

Notice that in this case the gauge group $\GG$ is equal to the
\emph{loop group} $\LL G$  of $G$, namely the group whose
elements are ordinary smooth maps $S^1\to G$ and whose operation
is pointwise multiplication. From this we see that since the space
$\Aa$ of connections is contractible then $B\Hh$ is weakly
homotopy equivalent to $B\LL G$. In any case we have the
following.

\begin{proposition} The groupoid $\Gg = (\Aa \times \LL G \rightrightarrows \Aa)$
is Morita equivalent to the groupoid $ \Hh = (G\times G
\rightrightarrows G)$
\end{proposition}

The Morita equivalence is obtained by taking the holonomy of the
connection. In view of this we will call $\Gg$ the \emph{inertia
groupoid} of the compact Lie group $G$. The groupoids $\Gg$, $\Hh$
and their Morita equivalence are relevant to the relation of the
twisted equivariant $K$-theory of $G$ (with respect to the
conjugate action) and the Verlinde algebra \cite{Freed,
LupercioUribe}.

\section{Smooth Manifolds}

\subsection{}
 Now we want to see that this
description when carried out on the groupoid associated to a
manifold $M$ matches the classical one of the loop space $\LL M$.
Let's denote by $\underline{M}$ and $\underline{S^1}$ the
groupoids associated to the trivial covers of $M$ and $S^1$
respectively ( i.e $\underline{M}_1=\underline{M}_0=M$ and
$\underline{S^1}_1=\underline{S^1}_0=S^1$) and by $\Mm$ the
groupoid associated to some fixed cover of $M$, i.e. if
$\{U_\alpha\}_\alpha$ is a cover of $M$ then
$$\Mm_1 := \bigsqcup_{(\alpha,\beta)} U_{\alpha\beta} \ \ \ \ \
\mbox{ and }  \ \ \ \ \ \ \Mm_0:= \bigsqcup_{\alpha} U_\alpha$$
with the natural source an target maps.\\
 The following
lemma will unveil the internal structure of $\Loop \Mm$.
\begin{lemma}
For $\Psi \in (\Loop \Mm)_0$ there exist a unique map $\psi:S^1
\to M$ such that the following diagram commutes
       $$
       \xymatrix{
       \Ss^1 \ar[r]^{P} \ar[d]_{\Psi}& \underline{S^1} \ar[d]^\psi \\
       \Mm \ar[r]_Q & \underline{M}
       }$$
where $P$ and $Q$ are the natural Morita morphisms induced by the
covers of $S^1$ and $M$ respectively.
\end{lemma}
\begin{proof}
There is no much choice in how to define $\psi$, we just need to
check that is well defined. Take a point $x$ in $S^1$, then every
pair of points in $(P_0)^{-1}\{x\}$ are joined by exactly one
arrow of $(\Ss^1)_1$, and all those arrows form precisely the set
$(P_1)^{-1}\{x\}$ (because $P$ is a Morita morphism). Then, under
the map $Q_1 \circ \Psi_1$, the set $(P_1)^{-1}\{x\}$ has to go
only one element $m$ in $\underline{M}_1$, this because the
groupoid $\underline{M}$ has only the identity arrows and no
relations among them. Hence $\psi(x)=m$ is well defined. The
continuity follows from the fact that $\Psi$ gives us the local
continuity of $\psi$.
\end{proof}
And so, we can get that
\begin{cor}
The induced map $\Loop Q :\Loop \Mm \to \Loop \underline{M}$
sends $\Psi$ to $\psi$, and for any other $\Psi'$ with $(\Loop
Q)_0(\Psi')=\psi$ then there exist $\Lambda$ in $(\Loop \Mm)_1$
joining $\Psi$ with $\Psi'$.
\end{cor}
\begin{proof}
We know from lemma \ref{moritaloopgroupoids} that $\Loop Q$ is a
Morita morphism; then as $(\Loop Q)_0(\Psi')=(\Loop
Q)_0(\Psi)=\psi$ and the fact that over $\psi$ there is the
identity arrow, there must exist an arrow $\Lambda$ in $(\Loop
\Mm)_1$ joining
 $\Psi'$ and $\Psi$.
\end{proof}

Let us abuse the notation and denote by $\LL M $ the groupoid
$\LL M \rightrightarrows \LL M$ with object space $\LL M$ and
only identity arrows.

\begin{proposition} The loop groupoid $\Loop \Mm$ is morita equivalent to $\LL M$.
\end{proposition}
\begin{proof}
By lemma \ref{moritaloopgroupoids} the map $\Loop Q :\Loop \Mm
\to \Loop \underline{M}$ is Morita. And by the previous lemmas is
clear that for any two admissible covers $W,W'' \in \WW$ the
groupoids $\Loop \underline{M}(W)$ and
 $\Loop \underline{M}(W')$ are equivalent. Hence for $\Loop \underline{M}$ it suffices to take the
trivial cover of $\real$ and then is clear that $\Loop
\underline{M}$ is equivalent to $\LL M$.
\end{proof}

We got, as expected, that the loop groupoid $\Loop \Mm$ carries
the same information as the loop space $\LL M$; this means that we
only need to consider the loops on $M$, in this case the
additional structure introduced by the groupoid representation of
the manifold doesn't give anything new.

\subsection{} This construction is compatible with the following fact. A
morphism of groupoids $\Psi :\Ss^1 \to \Mm$ induces a continuous
map on the classifying spaces
$$B \Ss^1 \to B \Mm$$
 where $B \Ss^1 \simeq S^1$ and $B \Mm \simeq M$ by a theorem of Segal \cite{Segal},
   thus producing a loop $\psi: S^1 \to M$.

\subsection{} Let's now briefly discuss the \emph{holonomy} of a gerbe with
connection as is done in the work  of Hitchin \cite{Hitchin}.
This will permit us define a line bundle with connection over the
loop space of $M$ and these ideas will be generalized for
groupoids.

For simplicity in what follows we are going to use subindices
instead of the groupoid notation. Let's recall what a gerbe with
connection means \cite[p.7]{Hitchin}: For a gerbe defined by a
cocycle $g_{\alpha\beta\gamma} : U_{\alpha\beta\gamma} \to \U{1}$,
a connection on it consist of a global 3-form $G \in
\Omega^3(M)$, 2-forms $F_\alpha \in \Omega^2(U_\alpha)$ and
1-forms $A_{\alpha\beta} \in \Omega^1(U_{\alpha\beta})$ such that
\begin{eqnarray*}
G|_{U_\alpha} & =& dF_\alpha\\
F_\beta - F_\alpha & = & dA_{\alpha\beta}\\
iA_{\alpha \beta} + iA_{\beta \gamma} + i A_{\gamma \alpha} & = &
g_{\alpha \beta \gamma}^{-1} dg_{\alpha \beta \gamma}
\end{eqnarray*}
 Adopting the definition of a gerbe as line bundles $\LLL_{\alpha\beta}$ on each
$U_{\alpha \beta}$ which in the triple intersection
$\LLL_{\alpha\beta} \LLL_{\beta\gamma} \LLL_{\gamma\alpha}$ get
trivialized by a cocycle $\theta_{\alpha\beta\gamma}$, we have
that a connection in this formalism is:
\begin{itemize}
\item a connection $\nabla_{\alpha\beta}$ on $\LLL_{\alpha \beta}$ such that
\item $\nabla_{\alpha\beta\gamma} \theta_{\alpha \beta \gamma} =0$
\item a 2-form $F_\alpha \in \Omega^2(U_\alpha)$ such that on $U_{\alpha \beta}$,
$F_\beta - F_\alpha= F_{\alpha\beta}$ is the curvature of
$\nabla_{\alpha \beta}$
\end{itemize}

We call the closed 3-form $G$ the \emph{curvature} of the gerbe
with connection, and we say that the connection gerbe is
\emph{flat} if the curvature vanishes. If this is the case,
$dF_\alpha =0$;  and using a Leray cover of $M$ (every
intersection of open sets is contractible) we can define 1-forms
$B_\alpha \in \Omega^1(U_\alpha)$ and 0-forms $f_{\alpha \beta}
\in \Omega^0(U_{\alpha \beta})$ such that:
$$F_\beta - F_\alpha = dA_{\alpha \beta} = d(B_\beta - B_\alpha)$$
$$A_{\alpha \beta} -B_\beta + B_\alpha = df_{\alpha \beta}$$
$$d(if_{\alpha \beta}  + if_{\beta \gamma} +if_{\gamma \alpha} - \log g_{\alpha \beta \gamma})=0$$
Then, as $\log g_{\alpha \beta \gamma}$ is only defined modulo
$2\pi i \integer$ we get a collection of constants $c_{\alpha
\beta \gamma}/2\pi \in \real / \integer$ that represents a
$\mathrm{\check{C}ech}$ class in $H^2(M,\real/\integer)$. This
2-cocycle is called the \emph{holonomy} of the connection.

If the holonomy is trivial ,then there are constants $k_{\alpha
\beta} \in 2 \pi i \integer$ such that $c_{\alpha \beta \gamma} =
k_{\alpha \beta} + k_{\beta \gamma} + k_{\gamma \alpha}$ and
defining $h_{\alpha \beta}:= \exp(if_{\alpha \beta} -k_{\alpha
\beta})$ we obtain a \emph{flat trivialization} of the gerbe
because $h_{\alpha\beta}h_{\beta\gamma}h_{\gamma\alpha}=g_{\alpha
\beta \gamma}$. If we have a second trivialization $h'_{\alpha
\beta}$, their difference $l_{\alpha \beta}:= h'_{\alpha \beta} /
h_{\alpha \beta}$ defines a line bundle $L$. Using that
$$iB_\beta - iB_\alpha -iA_{\alpha \beta} = d \log h_{\alpha \beta}$$
$$iB'_\beta - iB'_\alpha -iA_{\alpha \beta} = d \log h'_{\alpha \beta}$$
we deduce that
$$i(B' - B)_\beta - i(B' - B)_\alpha = d \log l_{\alpha \beta}$$
which defines a connection on $L$. As $d(B'_\alpha - B_\alpha)=0$
the curvature of this connection is zero; hence, the difference
of two flat trivializations of a gerbe is a flat line bundle.

Applying these ideas to a loop $f: S^1 \to M$  we get that, as
$S^1$ is 1-dimensional, the pull back of the gerbe with
connection to the circle is flat and has trivial holonomy. If we
identify flat trivializations that differ by a flat line bundle
with trivial holonomy, and we call it the moduli space of flat
trivializations, then to each loop,  this moduli space is being
acted freely and transitively by the moduli space of flat line
bundles $H^1(S^1,\real/\integer) \cong S^1$. Therefore we obtain
a principal $S^1$ bundle over $\LL M$.
 Now, if we take a path in the loop space
$$F:[0,1] \times S^1 \to M$$
the pull back of the gerbe connection is also flat and has
trivial holonomy, but it will give a canonical isomorphism
between the moduli space of flat trivializations of the gerbe on
$\{0\} \times S^1$ and $\{1\} \times S^1$. In this way we have
defined the parallel transport in the principal bundle, and hence
a connection on the line bundle over the loop space.

Thus from the curvature 3-form $G$ we obtained a 2-form on the
loop space, the curvature of the previous line bundle; this
process is known as \emph{transgression} and is what we are going
to generalize to orbifolds.

\section{Orbifolds}

\subsection{} We are going to analyze first the loop groupoid of a global quotient.
For $\Xx=[X/G]$ the groupoid given by the orbifold $X/G$ with $G$
a finite group and $X$ a connected $G$-space, we know from lemma
\ref{lemmaglobalquotient} that we only need to consider morphisms
$\Psi: \Ss^1 \to \Xx$ associated to the trivial cover of $\real$.
Now, from section \ref{morphismtrivialcover} we have that the
morphisms are determined by the image of $\real \times \{1\}$
under $\Psi_1$ because it lies on one connected component of
$\Xx_1$ (i.e. $\Psi_1(\real \times \{1\}) \subset X \times \{g\}$
for some $g \in G$). We know that $\Psi_0(x) = s(\Psi_1(x,1))$
and $\Psi_0(x+1) = t(\Psi_1(x,1))$, therefore if we consider the
pairs $(\psi,g)$ where $\psi:= \Psi_0$ and
$\Psi_1(\real\times\{1\}) \subset X \times \{g\}$ then

      $$
      \left\{ \begin{array}{c}
      \mbox{Morphisms } \Ss^1 \to \Xx \\
      \mbox{associated to the}\\
      \mbox{trivial cover.}
      \end{array} \right\} \stackrel{1-1}{\longleftrightarrow}
         \left\{ \begin{array}{c}
         \mbox{Pairs } (\psi,g) \mbox{ with} \\
         \psi:\real \to X \mbox{ and} \\
         \psi(x)\cdot g =\psi(x+1)
         \end{array} \right\} $$

An arrow $\Omega \in (\Loop \Xx)_1$ between $\Psi$ and $\Phi$ is
a map $\Omega: \real \times \integer \to \Xx_1$ where the
relevant information lies on $\real \times \{0\}$. As $X$ is
connected $\Omega(\real \times \{0\}) \subset X \times \{h\}$ for
some $h \in G$, and is easy to check that if $(\psi,g)$ and
$(\phi,k)$ are the pairs associated to $\Psi$ and $\Phi$
respectively, then
$$k =h^{-1} gh \ \ \ \ \ \mbox{and} \ \ \ \ \ \ \psi(x) \cdot h = \phi(x).$$
If we call $\PP_g$ the set of all pairs $(\psi,g)$
$$\PP_g = \{ (\psi,g) | \psi:\real \to X \ \ \& \ \ \psi(x)\cdot g = \psi(x+1)\}$$
then we can endow the set $\bigsqcup_g \PP_g$ with a natural $G$
action:
\begin{eqnarray*}
\bigsqcup_g \PP_g \times G & \to & \bigsqcup_g \PP_g\\
((\psi,g), h) & \mapsto & (\psi\cdot h, h^{-1}gh)
\end{eqnarray*}
so we can conclude
\begin{proposition}
If $\Xx=[X/G]$ is the quotient of the connected space $X$ by the
action of a finite group $G$, then the loop groupoid $\Loop \Xx$
is Morita equivalent to the groupoid
$$\LL [X/G]:=\begin{array}{c}
 \left( \bigsqcup_g \PP_g \right) \times G\\
\twodownarrows\\
\left( \bigsqcup_g  \PP_g \right)
\end{array}$$
given by the action explained before.
\end{proposition}
\begin{proof}
The functor $\Ff : \Loop [X/G] \to \LL [X/G]$ is given  by lemma
\ref{lemmaglobalquotient}. The Morita equivalence follows from
lemma \ref{arrowdeterminedbypoint} and the fact that  $\LL [X/G]$
is also a subgroupoid of $ \Loop [X/G]$ (with maps associated to
the trivial cover as in \ref{lemmaglobalquotient}) .
\end{proof}

By the same argument as in the case of the finite group $G$
acting on itself by conjugation, we obtain that
\begin{cor}
If $\Xx=[X/G]$ is the quotient of the connected space $X$ by the
action of a finite group $G$, then the loop groupoid $\Loop \Xx$
is Morita equivalent to the groupoid
$$\bigsqcup_{(g)} \left( \begin{array}{c}
   \PP_{g} \times C(g)\\
\twodownarrows\\
 \PP_g
\end{array} \right)$$
where $(g)$ runs over the conjugacy classes of elements in $G$.
\end{cor}

\subsection{} Now if we consider the groupoid $(\Loop \Xx)^\real$, the elements
$( \psi,g)$ in $\PP_g$ that are fixed under the action of $\real$
are the constant functions; so if $\psi(x)=x_0$ with $x_0$ a
fixed point in $X$, as
$$x_0 = \psi(x+1) =\psi(x) \cdot g = x_0 \cdot g$$
we have that $x_0$ belongs to the set $X^g$ of points in $X$
fixed by the action of $g$. Hence $(\PP_g)^\real \cong X^g$, and
using the previous corollary
\begin{proposition} \label{propositiontwistedsectorsglobalquotient}
The groupoid $(\Loop \Xx)^\real$ (or inertia groupoid $\Lambda
\Xx$) is Morita equivalent to the groupoid
$$\bigsqcup_{(g)} \left( \begin{array}{c}
   X^g \times \{g\}\times C(g)\\
\twodownarrows\\
 X^g \times \{g\}
\end{array} \right)$$
Therefore the groupoid $(\Loop \Xx)^\real$ is Morita equivalent
to the twisted sectors of the orbifold $X/G$.
\end{proposition}

\begin{rem} \label{remarkinertiagroupoid}
The inertia groupoid  $\wedge \Xx$ has the following
representations:
$$\bigsqcup_{(g)} \left( \begin{array}{c}
   X^g \times \{g\}\times C(g)\\
\twodownarrows\\
 X^g \times \{g\}
\end{array} \right)
\hookrightarrow
\begin{array}{c}
\left( \bigsqcup_g X^g \times \{g\}\right) \times G\\
\twodownarrows\\
 \bigsqcup_g X^g \times \{g\}
\end{array}
\stackrel{M}{\to} \bigsqcup_{(g)} \left( \begin{array}{c}
   X^g \times \{g\}\times C(g)\\
\twodownarrows\\
 X^g \times \{g\}
\end{array} \right)$$
where the first groupoid is Morita equivalent and a subgroupoid
of the second one.
\end{rem}

\subsection{} From \cite[Ex. 7.2.7]{LupercioUribe} we know that a gerbe $\LLL$ over
the groupoid $\overline{G}$ (with cocycle $\theta : G \times G
\to \U{1}$),
 as in section \ref{sectiongerbeongroup}, gives rise to a
gerbe $\LLL^\Xx$ over $\Xx$ in the natural way, i.e.
$$\begin{array}{cc}
\LLL^\Xx =&  X \times \LLL\\
 & \downarrow \\
\Xx_1 = & X \times G
\end{array}$$
with cocycle
\begin{eqnarray*}
\theta^\Xx : \Xx_2 = X \times G \times G & \to & \U{1}\\
(x,g,h) & \mapsto & \theta(g,h)
\end{eqnarray*}

Defining $\LLL^\Xx_g:=\LLL^\Xx |_{X\times \{g\}}$ is easy to see
that these line bundles are trivial, and using the notation of
section \ref{sectiongerbeongroup} we have that $\LLL^\Xx_g \cong
X \times \LLL_g$. Then is straightforward to generalize the
constructions and results of the section just mentioned to the
global quotient. If $\wedge \Xx$ is the inertia groupoid of $\Xx$
and $\rho$ is the map defined in lemma
\ref{lemmalinebundlediscretetorsion}, then

\begin{proposition}
The map
$$\rho^\Xx: \bigsqcup_{(g)} \left( X^g \times \{g\} \times C(g) \right)  \to  \U{1} $$
defined in $X^g \times \{g\} \times C(g)$ by
$$\rho^\Xx(x,g,h): = \rho(g,h)$$
is a morphism of groupoids $\rho^\Xx : \wedge \Xx \to
\overline{\U{1}}$. Therefore it determines a $C(g)$ equivariant
line bundle over $X^g$  for every $g$.
\end{proposition}

\begin{rem}
These line bundles precisely determine what Ruan \cite{Ruan}
called an inner local system of a discrete torsion $\theta \in
H^2(G, \U{1})$ over a global quotient $X/G$.
\end{rem}

\subsection{} We will now generalize the previous arguments to a
general \'{e}tale groupoid. For $\Gg$ an \'{e}tale groupoid, a
gerbe over $\Gg$ is determined by a 2-cocycle
$$\theta : \Gg_2= \Gg_1 \timests \Gg_1 \to \U{1}$$
which means that for $(a,b,c) \in \Gg_3 =\Gg_1 \timests \Gg_1
\timests \Gg_1$
$$\theta(a,b) \theta(ab,c) = \theta(a,bc)\theta(b,c);$$
this will give rise to a line bundle on the inertia groupoid.

An object $a$ in the inertia groupoid ($a \in \wedge \Gg_0$) is
an arrow of $\Gg_1$ such that its source and target are equal, so
$$\wedge\Gg_0 = \{ a \in \Gg_1 | s(a)=t(a)\}$$
and a morphism $v \in \wedge \Gg_1$ joining two objects $a,b$
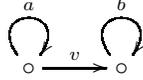
          $$\xymatrix{
          \circ \ar@(ul,ur)[]^a \ar[r]^v & \circ \ar@(ul,ur)[]^b}$$
is an arrow in $\Gg_1$ such that $a\cdot v = v\cdot b$; in some
way $b$ can be seen as the conjugate of $a$ by $v$ because $b =
v^{-1} \cdot a \cdot v$. But in order to identify the arrow $v$
we need to keep track of its source $a$, then we can consider the
morphisms of $\wedge \Gg$ as
$$\wedge \Gg_1 = \{ (a,v) \in \Gg_2 | a \in \wedge \Gg_0 \}$$
where $s(a,v) = a$ and $t(a,v) = v^{-1} \cdot a\cdot v$.

Let $\rho : \wedge \Gg_1 \to \U{1}$ be the map defined as
$$ (a,v) \stackrel{\rho}{\to} \frac{\theta(a,v)}{\theta(v,v^{-1}av)}$$
then this maps determines the mentioned line bundle.

\begin{lemma} \label{lemmamorphismofgroupoids}
The induced map $$\rho : \wedge \Gg \to \overline{\U{1}}$$ is a
morphism of groupoids.
\end{lemma}
\begin{proof}
Let $a,b,c \in \wedge \Gg_0$, and $v,w \in \wedge \Gg_1$
          $$
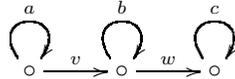
\xymatrix{
          \circ \ar@(ul,ur)[]^a \ar[r]^v & \circ \ar@(ul,ur)[]^b \ar[r]^w & \circ \ar@(ul,ur)[]^c}$$
such that $av = vb$ and $bw = wc$, we just need to prove that
$\rho(a,v) \rho (b,w) = \rho (a,vw)$. Using the cocycle condition
and noting that $b=v^{-1}av$ and $c=w^{-1}bw$,
 we get the following set of equalities:
\begin{eqnarray*}
\rho(a,vw)  = \frac{\theta(a,vw)}{\theta(vw,c)} & = &
 \frac{\theta(a,vw)}{\theta(av,w)} \frac{\theta(vb,w)}{\theta(v,bw)}
\frac{\theta(v,wc)}{\theta(vw,c)}\\
& = &\frac{\theta(a,v)}{\theta(v,w)}
\frac{\theta(b,w)}{\theta(v,b)}
\frac{\theta(v,w)}{\theta(w,c)}\\
& = & \frac{\theta(a,v)}{\theta(v,b)}
\frac{\theta(b,w)}{\theta(w,c)} = \rho(a,v) \rho(b,w)
\end{eqnarray*}
The  continuity is clear.
\end{proof}
So we get an interesting result:
\begin{theorem}
Let $\Gg$ be an \'{e}tale groupoid and $\wedge \Gg$ its inertia
groupoid. A gerbe over $\Gg$ determines a morphism of groupoids
$\rho :\wedge \Gg \to \overline{\U{1}}$; hence
 a line bundle over the inertia groupoid $\wedge \Gg$.
\end{theorem}

When the groupoid $\Gg$ is an orbifold its is known that the
twisted sectors $\widetilde{\Sigma_1 \Gg}$ is a subgroupoid of
$\wedge \Gg$ which is at the same time Morita equivalent to it
(see Appendix).
 Then the morphism $\rho$ induces a line orbibundle over $\widetilde{\Sigma_1 \Gg}$ in the
natural way.

\begin{cor}
A gerbe over an orbifold $\Gg$ (in the sense of
\cite{LupercioUribe}) determines an inner local system over the
twisted sectors $\widetilde{\Sigma_1 \Gg}$.
\end{cor}
 \begin{proof}
Calling $\eta : \wedge \Gg_0 \to \Gg_1$ the inclusion, and $\LL$
the line bundle over $\Gg_1$ determined by the gerbe, the line
bundle on $\wedge \Gg$ is constructed from $\eta^*\LL$ and the
attaching maps given by $\rho$. The properties of the gerbe
\cite{LupercioUribe}:
\begin{itemize}
\item $i^* \LL \cong \LL^{-1}$
\item $\pi_1^* \LL \cdot \pi_2^* \LL \cdot (im)^* \LL \cong 1$
\end{itemize}
clearly imply the conditions in \cite[Def. 3.1.6]{Ruan} for an
inner local system.
\end{proof}

\section{Appendix}
\subsection{}In this section we will summarize some facts regarding
orbifolds and its twisted sectors, and the relation with the
groupoid approach. For a more detailed account on orbifolds we
recommend \cite{LupercioUribe, Ruan}; the notation and results
that will be used in what follows can be seen in \cite[Sect.
5.1]{LupercioUribe}.

\subsection{}Let $X$ be an orbifold and $\{(V_p, G_p, \pi_p)\}_{p \in
X}$ its orbifold structure, the groupoid $\Xx$ associated to it
(cf. \cite{MoerdijkPronk1}) consist of the following information:
\begin{itemize}
\item Objects ($\Xx_0$): $\bigsqcup_{p \in X} V_p$.
\item Morphisms ($\Xx_1$): an arrow $r : (x_1,V_1) \to (x_2,V_2)$ is an equivalence class
of triples $r= [\lambda_1, w , \lambda_2]$ where $w \in W$ for
another orbifold chart $(W,H,\rho)$, and the $\lambda_i$'s are
embeddings $\lambda_i : W \to V_i$ which are $\phi_i$
equivariant, with $\phi_i: H \to G_i$ monomorphisms;
          $$\xymatrix{
            (V_1,G_1) & & (V_2,G_2)\\
            & (W,H) \ar[lu]^{(\lambda_1, \phi_1)} \ar[ru]_{(\lambda_2, \phi_2)} &
            }$$
\end{itemize}

\subsection{}If  $r \in \wedge \Xx_0$ is an element in the inertia
groupoid, then its source and target must be equal. In particular
we would have that $x_1 =x_2(=x)$ and that $V_1 = V_2(=V)$. Since
the only automorphisms of the orbifold chart $(V,G,\pi)$ come from
conjugation, then there exist $g \in G$ such that $\lambda_2 =
\lambda_1 \cdot g$ and $\phi_2 = g^{-1}\cdot  \phi_1 \cdot g$.
Now it's easy to see that
$$[\lambda_1, w, \lambda_2 ]=[id_V, x, \lambda_g]$$
where $\lambda_g : V \to V$, $\lambda_g(y)=y \cdot g$ and $\phi_g
: G \to G$ is the map obtained conjugating by $g$. So $x \in V^g$
is a fixed point under the action of $g$ in $V$. In this way we
can see that the objects of the inertia groupoid can be seen as
the fixed point sets of the actions of the elements of the group
$G_p$ in $V_p$
$$\Xx_0 \cong \bigsqcup_{p \in X} \left( \bigsqcup_{g \in G_p} (V_p)^g \times \{g\}\right)$$
Again, using that the only automorphisms of charts are given by
conjugation, we can see that if we restrict our attention to the
morphisms that are defined on
$$Z_p:= \bigsqcup_{g \in G_p} (V_p)^g \times \{g\}$$
we get the same picture as in proposition
\ref{propositiontwistedsectorsglobalquotient}. Hence, the inertia
groupoid $\wedge \Xx$ restricted $Z_p$ is Morita equivalent to
$$ \wedge \Xx |_{Z_p} \stackrel{M}{\cong} \bigsqcup_{(g) \subset G_p}
\left( \begin{array}{c}
(V_p)^g \times \{g\} \times C_{G_p}(g)\\
\twodownarrows\\
(V_p)^g \times \{g\}
\end{array} \right)$$
and this is the description given in \cite{Ruan} of the twisted
sectors $\widetilde{\Sigma_1 X}$. Therefore we can conclude that
\begin{proposition}
Let $X$ be an orbifold whose associated groupoid is $\Xx$. Then
the inertia groupoid $\wedge \Xx$ is Morita equivalent to the
twisted sectors $\widetilde{\Sigma_1 X}$. Moreover,
$\widetilde{\Sigma_1 X}$ is also a subgroupoid of $\wedge \Xx$:
      $$\xymatrix{
      \widetilde{\Sigma_1 X} \ar@{^{(}->}[r] & \wedge \Xx \ar[r]^M & \widetilde{\Sigma_1 X}}$$
\end{proposition}
\begin{proof}
It follows directly from the previous discussion and the case of a
global quotient (remark \ref{remarkinertiagroupoid}).
\end{proof}

\bibliographystyle{amsplain}
\bibliography{loopgroup}

\providecommand{\bysame}{\leavevmode\hbox to3em{\hrulefill}\thinspace}
\providecommand{\MR}{\relax\ifhmode\unskip\space\fi MR }
\providecommand{\MRhref}[2]{%
  \href{http://www.ams.org/mathscinet-getitem?mr=#1}{#2}
}
\providecommand{\href}[2]{#2}
\begin{thebibliography}{10}

\bibitem{Adem}
A.~Adem, \emph{Characters and k-theory of discrete groups}, Invent. Math.
  \textbf{114} (1993), 489--514.

\bibitem{BridsonHaefliger}
Martin~R. Bridson and Andr{\'e} Haefliger, \emph{Metric spaces of non-positive
  curvature}, Springer-Verlag, Berlin, 1999.

\bibitem{Chen}
W.~Chen, \emph{A homotopy theory of orbispaces}, arXivmath.AT/0102020: (2001).

\bibitem{CrainicMoerdijk}
Marius Crainic and Ieke Moerdijk, \emph{A homology theory for \'etale
  groupoids}, J. Reine Angew. Math. \textbf{521} (2000), 25--46.

\bibitem{Freed}
D.~Freed, \emph{{The Verlinde algebra is twisted equivariant K-theory}},
  arXiv:math.RT/0101038.

\bibitem{Hitchin}
N.~Hitchin, \emph{Lectures notes on special lagrangians submanifolds},
  arXiv:math.DG/9907034 (1999).

\bibitem{Kawasaki}
T.~Kawasaki, \emph{The signature theorem for v-manifolds}, Topology \textbf{17}
  (1978), 75--83.

\bibitem{LupercioUribe}
E.~Lupercio and B.~Uribe, \emph{Gerbes over orbifolds and {K}-theory},
  arXiv:math.AT/0105039 (2001).

\bibitem{Moerdijk1}
I.~Moerdijk, \emph{Calssifying topos and foliations}, Ann. Inst. Fourier
  \textbf{41} (1991), 189--209.

\bibitem{Moerdijk2}
\bysame, \emph{Orbifolds as groupoids: an introduction.}, arXiv:math.DG/0203100
  (2002).

\bibitem{MoerdijkPronk1}
I.~Moerdijk and D.~Pronk, \emph{Orbifolds, sheaves and grupoids}, K-Theory
  \textbf{12} (1997), 3--21.

\bibitem{Mrcun}
Janez Mr{\v{c}}un, \emph{An extension of the {R}eeb stability theorem},
  Topology Appl. \textbf{70} (1996), no.~1, 25--55.

\bibitem{Pronk}
D.~Pronk, \emph{Etendues and stacks as bicategories of fractions}, Comm. Math.
  \textbf{102} (1996), 243--303.

\bibitem{Ruan}
Y.~Ruan, \emph{Stringy geometry and topology of orbifolds},
  arXiv:math.AG/0011149 (2000).

\bibitem{Segal}
G.~Segal, \emph{Classifying spaces and spectral sequences}, Inst. Hautes
  {E}tudes Sci. Publ. Math. \textbf{34} (1968), 105--112.

\end{thebibliography}

\end{document}